\documentclass[12pt,twoside]{amsart}
\usepackage{amssymb}
\usepackage{graphicx}

\nonstopmode

\numberwithin{equation}{section}
\font\tengothic=eufm10 scaled\magstep 1
\font\sevengothic=eufm7 scaled\magstep 1
\newfam\gothicfam
          \textfont\gothicfam=\tengothic
          \scriptfont\gothicfam=\sevengothic

\DeclareMathOperator{\chara}{char}

\def\cocoa{{\hbox{\rm C\kern-.13em o\kern-.07em C\kern-.13em o\kern-.15em A}}}

\DeclareMathOperator{\pnt}{\raise 0.5mm \hbox{\large\bf.}}

\newtheorem{theorem}{Theorem}[section]
\newtheorem{lemma}[theorem]{Lemma}
\newtheorem{proposition}[theorem]{Proposition}
\newtheorem{corollary}[theorem]{Corollary}
\newtheorem{conjecture}[theorem]{Conjecture}

\newtheorem{question}[theorem]{Question}
\newtheorem{problem}[theorem]{Problem}

\theoremstyle{definition}
\newtheorem{definition}[theorem]{Definition} 
\newtheorem{remark}[theorem]{Remark}
\newtheorem{example}[theorem]{Example}

\begin{document}

\title{A tour of the Weak and Strong Lefschetz Properties}

\author[Juan Migliore]{Juan Migliore}
\address{
Department of Mathematics, University of Notre Dame, Notre Dame, IN
46556, USA}
\email{Juan.C.Migliore.1@nd.edu}
\author[Uwe Nagel]{Uwe Nagel}
\address{Department of Mathematics,
University of Kentucky, 715 Patterson Office Tower,
Lexington, KY 40506-0027, USA}
\email{uwe.nagel@uky.edu}

\thanks{
This work was partially supported by two grants from the
Simons Foundation (\#208579 to Juan Migliore and \#208869 to Uwe Nagel).
}

\begin{abstract}
An artinian graded algebra, $A$, is said to have the Weak Lefschetz property (WLP) if multiplication by a general linear form has maximal rank in every degree.  A vast quantity of work has been done studying and applying this property, touching on numerous and diverse areas of algebraic geometry, commutative algebra, and combinatorics.  Amazingly, though, much of this work has a ``common ancestor" in a theorem originally due to Stanley, although subsequently reproved by others.  In this paper we describe the different directions in which research has moved starting with this theorem, and we discuss some of the open questions that continue to motivate current research.

\end{abstract}

\maketitle

\tableofcontents

\section{Introduction}

The Weak and Strong Lefschetz properties are strongly connected to
many topics in algebraic geometry, commutative algebra and
combinatorics.  Some of these connections are quite surprising and
still not completely understood, and much work remains to be done.
In this expository paper we give an overview of known results on the
Weak and Strong Lefschetz properties, with an emphasis on the vast
number of different approaches and tools that have been used, and
connections that have been made with seemingly unrelated problems.
One goal of this paper is to illustrate the variety of methods  and
connections that have been brought to bear on this problem for
different families of algebras.  We also discuss open problems.

Considering the amazing breadth and depth of the results that  have
been found on this topic, and the tools and connections that have
been associated with it, it is very interesting to note that to a
large degree, one result motivated this entire area.  This result is
the following.  It was proved by R. Stanley \cite{stanley} in 1980 using algebraic topology, by J.\ Watanabe in 1987 using representation theory, by Reid, Roberts and Roitman \cite{RRR} in 1991 with algebraic methods, by Herzog and Popescu \cite{HP} (unpublished) in 2005 essentially with linear algebra, and it follows from work of Ikeda \cite{sekiguchi} in 1996 using combinatorial methods.

\begin{theorem}
\label{swrrr thm}

Let $R = k[x_1,\dots,x_r]$, where $k$ has  {\em characteristic zero}.
Let $I$ be an artinian monomial complete intersection, i.e.
\[
I = \langle x_1^{a_1},\dots,x_r^{a_r} \rangle.
\]
Let $\ell$ be a general linear form.  Then for any positive
integers $d$ and $i$, the homomorphism induced by multiplication by
$\ell^d$,
\[
\times \ell^d : [R/I]_i \rightarrow [R/I]_{i+d}
\]
has maximal rank.  (In particular, this is true when $d=1$.)
\end{theorem}

This paper is organized around the ways that subsequent research
owes its roots to this theorem.

Our account is by no means exhaustive. Fortunately, the manuscript \cite{HMMNWW} has appeared recently. It gives an overview of the Lefschetz properties from a different perspective, focusing more on the local case, representation theory, and combinatorial connections different from those presented here.

There is one topic that is neither treated in \cite{HMMNWW} nor here but that is worth mentioning briefly.  
 In \cite{MMN}, examples of monomial ideals  were exhibited that did not have the WLP,  but that could be deformed to ideals with the WLP. A systematic way for producing such deformations that preserve the Hilbert function has been proposed by Cook and the second author in \cite{CN-JPPA}. The idea is to lift the given monomial ideal to the homogenous ideal of a set of points and then pass to a general hyperplane section of the latter. It is shown in \cite{CN-JPPA} that this procedure does indeed produce ideals with the WLP for a certain class of monomial ideals without the WLP.
\smallskip

In May 2011, the first author gave a talk at the Midwest Commutative Algebra \& Geometry conference at Purdue University on this topic.  This paper is a vast expansion and extension of that talk, containing many more details and several new topics.  The first author is grateful to the organizers of the conference for the invitation to speak, and both authors are grateful to Hal Schenck for his suggestion that they write this paper.

The authors also thank David Cook II for helpful comments.


\section{Definitions and background}

Let $k$ be an infinite field.    We will often take char($k) = 0$,
but we will see that changing the characteristic produces
interesting new questions  (and even more interesting answers!).

Let $R = k[x_1,\dots,x_r]$  be the graded polynomial ring in $r$
variables over $k$. Let
\[
A = R/I = \bigoplus_{i=0}^n A_i
\]
be a graded artinian algebra.  Note that $A$ is finite dimensional
over $k$.

\begin{definition}
For any standard graded algebra $A$ (not necessarily artinian), the {\em Hilbert function} of $A$ is the function
\[
\underline{h}_A : \mathbb N \rightarrow \mathbb N
\]
defined by $\underline{h}_A(t) = \dim [A]_t$.  One can express $\underline h_A$ as a sequence
\[
(h_0 = 1, h_1, h_2, h_3, \dots).
\]
An {\em $O$-sequence} is a sequence of positive integers that occurs as the Hilbert function of some graded algebra.  When $A$ is Cohen-Macaulay, its {\em $h$-vector} is the Hilbert function of an artinian reduction of $A$.  In particular, when $A$ is artinian, its Hilbert function is equal to its $h$-vector.
\end{definition}

\begin{definition}
An {\em almost complete intersection} is a standard graded algebra $A = R/I$ which is Cohen-Macaulay, and for which the number of minimal generators of $I$ is one more than its codimension.
\end{definition}

\begin{definition}
$A$ is {\em level of Cohen-Macaulay type $t$} if its socle is
concentrated in one degree (e.g. a complete intersection) and has
dimension $t$.
\end{definition}

\begin{definition}  \label{def of wlp slp}
Let $\ell$ be a general linear form.    We say that $A$ has  the
{\em Weak Lefschetz Property (WLP)} if the homomorphism induced by
multiplication by $\ell$,
\[
\times \ell : A_i \rightarrow A_{i+1},
\]
 has maximal rank for all $i$ (i.e. is  injective or surjective).
We say that $A$ has the  {\em Strong Lefschetz Property (SLP)} if
\[
\times \ell^d : A_i \rightarrow A_{i+d}
\]
has maximal rank for all $i$ and $d$ (i.e. is  injective or
surjective).
\end{definition}

\begin{remark}
\begin{itemize}
\item[(a)] One motivation for the work described in this paper is
that something interesting should be going on if multiplication
by  a general linear form does not induce a homomorphism of
maximal rank, even in one degree.

\item[(b)] Later we will see that there is a strong  connection
to Fr\"oberg's conjecture.  In this regard, we note that
$\ell^d$  should not be considered to be a ``general" form of
degree $d$, since in the vector space $[R]_d$ ($d>1$), those
forms that are pure powers of linear forms form a proper
Zariski-closed subset.

\item[(c)]  Suppose that $\deg f = d$ and
$\times f : [R/I]_i \rightarrow [R/I]_{i+d}$ has maximal rank,
for  all $i$.   Pardue and Richert \cite{PR} call such an $f$
{\em semi-regular}.   Reid, Roberts and Roitman \cite{RRR} call
such an $f$ {\em faithful}. If $\times f^j : [R/I]_i \rightarrow
[R/I]_{i+dj}$ has maximal rank for all $i$ and all $j$, they
call such an $f$ {\em strongly faithful}. So $R/I$ has the WLP
if $R$ contains a linear faithful element, and $R/I$ has the SLP
if $R$ contains a linear strongly faithful element.

\item[(d)] Several authors consider the question of the ranks that
arise if $\times \ell^d$ is replaced by $\times F$ for a general
$F$  of degree $d$.  This is the essence of the Fr\"oberg
conjecture, is related to the WLP, and will be discussed below
in section \ref{froeberg section}.

\end{itemize}
\end{remark}

\bigskip

How do we determine if $R/I$ fails to have the WLP?   Let $\ell$ be
a general linear form and fix an integer $i$.    Then we have an
exact sequence
\[
[R/I]_{i-1} \stackrel{\times \ell}{\longrightarrow} [R/I]_i \rightarrow [R/(I,\ell)]_i \rightarrow 0.
\]
Thus $\times \ell$ fails to have maximal rank from degree $i-1$ to
degree $i$ if and only if
 \[
 \dim [R/(I,\ell)]_i > \max \{ 0, \dim[R/I]_i - \dim [R/I]_{i-1} \}.
 \]
More precisely, if we want to show that the WLP fails, it is enough
to identify a degree $i$ for which we can produce one of the
following two pieces of information:

\medskip

\begin{itemize}
\item[(i)] $\dim [R/I]_{i-1} \leq \dim [R/I]_{i}$ and $\dim  [R/(I,\ell)]_i >  \dim[R/I]_i - \dim [R/I]_{i-1}$; in this case we loosely say that WLP {\em fails because of injectivity}; or \\

\item[(ii)] $\dim [R/I]_{i-1} \geq \dim [R/I]_{i}$ and $\dim  [R/(I,\ell)]_i > 0$; in this case we loosely say that WLP {\em fails because of surjectivity}.
\end{itemize}

\noindent In general, even  identifying which $i$ is the correct place to look can be difficult. Then  determining which of (i) or (ii) holds, and establishing both inequalities, is often very challenging. This is where computer algebra programs have been very useful, in suggesting where to look and what to look for!  On the other hand, to prove that $R/I$ {\em does} have the WLP, the following result is helpful:

\begin{proposition}[\cite{MMN}, Proposition 2.1]  \label{MMN tool}
Let $R/I$ be an artinian standard graded algebra and let $\ell$ be a general linear form.  Consider the homomorphisms $\phi_d : [R/I]_d \rightarrow [R/I]_{d+1}$ defined by multiplication by $\ell$, for $d \geq 0$.

\begin{itemize}
\item[(a)] If $\phi_{d_0}$ is surjective for some $d_0$ then $\phi_d$ is surjective for all $d \geq d_0$.

\item[(b)] If $R/I$ is level and $\phi_{d_0}$ is injective for some $d_0 \geq 0$ then $\phi_d$ is injective for all $d \leq d_0$.

\item[(c)] In particular, if $R/I$ is level and $\dim [R/I]_{d_0}
= \dim [R/I]_{d_0+1}$ for some $d_0$ then $R/I$ has the WLP if and
only if $\phi_{d_0}$ is injective (and hence is an isomorphism).
\end{itemize}
\end{proposition}

\noindent  This result helps to narrow down where one has to look, especially in the situation where we want to show that the WLP {\em does} hold.    In this case  you have to find the critical degrees and {\em then} show that surjectivity {\em and} (usually) injectivity {\em do} hold just on two (or occasionally one) spot.

In the case of one variable, the WLP and SLP are trivial since all ideals are principal.  The case of two variables also has a nice result, at least in characteristic 0:

\begin{theorem}[\cite{HMNW}] \label{hmnw 2 var}

If $\hbox{\em char }(k) = 0$ and $I$ is {\em any} homogeneous ideal
in $k[x,y]$ then $R/I$ has the SLP.
\end{theorem}

The proof of this result used generic initial ideals with respect to
the reverse lexicographic order.   In the case of the WLP, it is not
hard to show that the above theorem is true  in any characteristic
(\cite{MZ}, \cite{LZ}, \cite{CN2}). However, the characteristic zero
assumption cannot be omitted for guaranteeing the SLP. In fact, also
the WLP may fail if there are at least three variables. The
following is an easy exercise:

\begin{lemma}
  \label{lem:ci-pos-char}
Assume $\hbox{char}(k) = p$.  Consider the ideal
\[
I = \langle x_1^p ,\dots,x_r^p \rangle \subset R = k[x_1,\dots,x_r],
\]
where  $r \geq 2$. Then

\begin{itemize}

\item $R/\langle x_1^p,\dots,x_r^p \rangle$ fails the SLP for all
 $r \geq 2$ (except if $r=p=2$, where  $R/\langle
 x^4,y^4\rangle$ fails the SLP).

\vspace{.1in}

\item It fails the WLP for all $r \geq 3$.

\vspace{.1in}

\item It {\em has} the WLP when $r=2$.

\end{itemize}

\vspace{.2in}
\end{lemma}

In Section \ref{pos char} we will discuss the presence of the WLP in
positive characteristic in more detail.

A useful consequence of knowing that an algebra $A$ has the WLP or
SLP is that its Hilbert function is {\em unimodal}.  In fact,  the
Hilbert functions of algebras with the WLP have been completely
classified:

\begin{proposition}[\cite{HMNW}]
Let $\underline{h} = (1,h_1,h_2,\dots,h_s)$ be a finite sequence of positive integers.  Then $\underline{h}$ is the Hilbert function of a graded artinian algebra with the WLP if and only if the positive part of the first difference is an $O$-sequence and after that the first difference is non-positive until $\underline{h}$ reaches 0.  Furthermore, this is also a necessary and sufficient condition for $\underline{h}$ to be the Hilbert function of a graded artinian algebra with the SLP.
\end{proposition}

The challenge is thus to study the WLP and SLP (or their failure),
and the behavior of the Hilbert function, for  interesting {\em
families} of algebras.  Most of the results below fall into this
description. It should also be noted that conversely, some Hilbert
functions  $\underline{h}$ {\em force} any algebra with Hilbert
function $\underline h$ to have the WLP -- these were classified in
\cite{MZ}.

In the rest of this paper, we indicate different directions of
research that have been motivated by Theorem \ref{swrrr thm}; in
most cases, there also remain many intriguing open problems.


\section{Complete intersections and Gorenstein algebras}

By semicontinuity, a consequence of Theorem \ref{swrrr thm} is  that
a {\em general} complete intersection with fixed generator degrees
has the WLP and the SLP.

\begin{question} \label{ci question}
Do {\em all} artinian complete intersections have the WLP or the
SLP in characteristic 0?
\end{question}

We know that the answer is trivially ``yes" in one and two
variables.   In three or more variables, the following is the most
complete result known to date.

\begin{theorem}[\cite{HMNW}]
Let $R = k[x,y,z]$, where $\chara (k) = 0$.  Let $I = \langle F_1,F_2,F_3 \rangle$
be a complete intersection.  Then   $R/I$ has the WLP.
\end{theorem}

The proof  of this result introduced the use of the syzygy module of
$I$, and its sheafification, the syzygy bundle.  Subsequently,
several papers have used the syzygy module to study the WLP for
different kinds of ideals (see, e.g., \cite{BK2-char p}, \cite{BK-WLP},
\cite{ss}, \cite{HSS}, \cite{MMN}, \cite{CN2}). In
the  case of complete intersections in $k[x,y,z]$, the syzygy bundle has rank 2.    The WLP is almost
immediate in the ``easy" cases, and semistability and the
Grauert-M\"ulich theorem give the needed information about
$R/(I,\ell)$ in the ``interesting" cases.

\begin{remark}
\begin{enumerate}

\item[(i)] The SLP is still wide open for complete intersections
in three or more variables, and in fact even the WLP is open for complete intersections
of  arbitrary codimension $\geq 4$.  Some partial results on the WLP for arbitrary complete intersections in four variables have been obtained recently by the authors together with Boij and Mir\'o-Roig, in work in progress.

\item[(ii)] It was conjectured by Reid, Roberts and Roitman \cite{RRR} that the
answer to both parts of Question \ref{ci question} is yes.

\end{enumerate}
\end{remark}

We have seen that conjecturally (and known in special cases), all
complete intersections have the WLP.   Complete intersections are a
special case of Gorenstein algebras.  Does the conjecture extend to
the Gorenstein case?  That is,

\begin{question}\label{gor wlp q} Do all graded artinian Gorenstein
algebras have the WLP?  If not, what are classes of artinian
Gorenstein algebras that do have this property?

\end{question}

The answer to the first question is a resounding ``no."  Indeed,
Stanley \cite{stanley-adv} in 1978 gave an example of an artinian
Gorenstein algebra with Hilbert function $(1, 13, 12, 13, 1)$, which
because of the non-unimodality clearly does not have the WLP.  Other
examples of non-unimodality for Gorenstein algebras were given by
Bernstein and Iarrobino \cite{BI}, by Boij \cite{Bo} and by Boij and
Laksov \cite{BL}. Even among Gorenstein algebras with unimodal
Hilbert functions, WLP  does not necessarily hold.  For instance, an
example in codimension 4 was given by Ikeda \cite{ikeda} in 1996.

On the other hand, the problem  in three variables is still wide
open,  with only special cases known (see for instance \cite{MZ2},
\cite{amasaki}):

\begin{question}
Does {\em every} artinian Gorenstein quotient of $k[x,y,z]$ have the
WLP, provided $\chara (k) = 0$?  What about the SLP?
\end{question}

\noindent Given the complete intersection result for three variables
mentioned above, this is a very natural and intriguing question.

In four variables, the result of Ikeda mentioned above shows that
WLP need not hold.  Nevertheless, the main result of \cite{MNZ}
shows that for small initial degree, the Hilbert functions are still
precisely those of Gorenstein algebras with the WLP.  More
precisely, it was shown that if the $h$-vector is
$(1,4,h_2,h_3,h_4,\dots)$ and $h_4 \leq 33$ then this result holds.
More recently, using the same methods, Seo and Srinivasan \cite{ss} extended this to $h_4 = 34$.   Thus, the
result holds for initial degree $\leq 4$.
\smallskip

Another interesting special case is the situation in which the generators of the ideal have small degree. We say that an algebra $R/I$ is presented by quadrics if the ideal $I$ is generated by quadrics. Such ideals occur naturally, for example, as homogeneous ideals of sufficiently positive embeddings of smooth
projective varieties (\cite{EL}) or as Stanley-Reisner ideals of
simplicial flag complexes (\cite{St-book}). Gorenstein algebras presented by quadrics are studied, for example,  in \cite{MN-gor-quad}. There, the following conjecture has been proposed.

\begin{conjecture}[\cite{MN-gor-quad}]
  \label{conj:Gor-quadrics}
Any artinian Gorenstein algebra presented by quadrics, over a field $k$ of characteristic zero, has the WLP.
\end{conjecture}

The conjecture predicts in particular that if the socle degree is at least 3 then the multiplication by a general linear form from degree one to degree two is injective.  Though this is established in some cases in \cite{MN-gor-quad}, even this special case of the conjecture is open.

The analog of Question \ref{gor wlp q} is also of interest for rings of positive dimension. If $A$ is a Gorenstein ring of dimension $d$, then $A$ is said to have the WLP if a general artinian reduction of $A$ has the WLP, that is, if $A/\langle L_1,\ldots,L_d \rangle$ has the WLP, where $L_1,\ldots,L_d \in A$ are general forms of degree 1. Recall that the Stanley-Reisner ring of the boundary complex of a convex polytope is a reduced Gorenstein ring. The so-called g-theorem classifies their Hilbert functions. The necessity of the conditions on the Hilbert function is a consequence of the following result by Stanley.

\begin{theorem}[\cite{St-faces}]
  \label{thm:g-thm}
The Stanley-Reisner ring of the boundary complex of a convex polytope over a field $k$ has the SLP if $\chara (k) = 0$.
\end{theorem}

The so-called g-conjecture states that the mentioned conditions on the Hilbert function characterize in fact the Hilbert functions of the Stanley-Reisner rings of triangulations of spheres. Note that there are many more such triangulations than boundary complexes of convex polytopes. In this regard, the following question  merits highlighting:

\begin{question}
  \label{q:WLP-points}
Does  a {\em reduced}, arithmetically
Gorenstein set of points in $\mathbb P^n$ have the WLP, provided $\chara (k) = 0$?
\end{question}

We point out that if this question has an affirmative answer, then,
by  the main result of \cite{MN3}, we have a classification of the
Hilbert functions of reduced, arithmetically Gorenstein schemes:
their $h$-vectors are precisely the SI-sequences, meaning that they
are symmetric, with the first half itself a differentiable
$O$-sequence.

An affirmative answer to Question \ref{q:WLP-points} would  also
imply the g-conjecture, thus providing a characterization of the
face vectors of triangulations of a sphere. Moreover, the methods
used to establish the WLP could lead to information about the face
vectors of triangulations of other manifolds as well. In fact, Novik
and Swartz \cite{NS-Bu}, Theorem 1.4, show that a certain quotient
of the Stanley-Reisner ring of any orientable $k$-homology manifold
without boundary is a Gorenstein ring. Kalai conjectured that this
Gorenstein ring has the SLP. If true, this would establish new
restrictions on the face vectors of these complexes. A special case
of Kalai's conjecture has been proven in Theorem 1.6 of
\cite{NS-Bu}.


\section{Monomial level algebras}

Note that $R/ \langle x_1^{a_1},\dots,x_r^{a_r} \rangle$ is also a
level artinian monomial algebra.

\begin{question}
  \label{q:failure-WLP-level}
Which (if any) level artinian monomial algebras fail the WLP or SLP?
\end{question}

The first result in this direction is a positive one:

\begin{theorem}[Hausel \cite{Ha}, Theorem 6.2] \label{hausel}
Let $A$ be a monomial artinian level algebra of socle degree $e$. If the field $k$ has characteristic zero, then for a {\em general} linear form $L$, the induced multiplication
\[
\times L : A_j \rightarrow A_{j+1}
\]
is an injection, for all $j = 0,1,\dots,\lfloor \frac{e-1}{2} \rfloor$.  In particular, over any field the sequence
\[
1, h_1 -1, h_2 - h_1 ,\dots, h_{\lfloor \frac{e-1}{2} \rfloor +1} - h_{\lfloor \frac{e-1}{2} \rfloor}
\]
is an $O$-sequence, i.e.\ the ``first half'' of $\underline{h}$ is a {\em differentiable $O$-sequence}.
\end{theorem}

Thus roughly ``half" the algebra does satisfy the WLP.  What about the second half?
The first counterexample was due to Zanello (\cite{Z-nonunimodal} Example
7),  who showed that the WLP does not necessarily hold for monomial
level algebras even in three variables.  His example had $h$-vector
$(1,3,5,5)$.  Subsequently, Brenner and Kaid (\cite{BK-WLP} Example
3.1) produced an example of a level artinian monomial {\em almost
complete intersection} algebra that fails the WLP; this algebra has
$h$-vector $(1,3,6,6,3)$ and in particular, Cohen-Macaulay type 3.
The study of such almost complete intersections was continued by
Migliore, Mir\'o-Roig and Nagel \cite{MMN}, and more recently by
Cook and Nagel \cite{CN}, \cite{CN2} (see also Section \ref{pos char}).

The Hilbert functions of the algebras considered in Question \ref{q:failure-WLP-level} are of great interest in a number of areas. In fact, they are better known under a different name.

\begin{definition}
A {\em pure $O$-sequence of type $t$ in $r$ variables} is the
Hilbert function of a level artinian monomial algebra
$k[x_1,\dots,x_r]/I$ of Cohen-Macaulay type $t$.
 \end{definition}

\begin{question}
We have already seen that level artinian monomial algebras do not
necessarily have the WLP.  Nevertheless, are their Hilbert functions
unimodal?  That is, are all pure $O$-sequences unimodal?   If not,
can we find subfamilies, depending on the type $t$ and/or the number
of variables $r$, that are unimodal?  And if they are not
necessarily unimodal, ``how non-unimodal" can they be?
\end{question}

\begin{remark}
If $I$ is a {\em monomial ideal} in $R = k[x_1,\dots,x_r]$  then the linear form $\ell = x_1+\cdots +x_r$ is ``general enough" to determine if $R/I$ has the WLP or SLP. This observation has been extremely useful in simplifying calculations to show the existence or failure of the WLP.
In   \cite{MMN}, Proposition 2.2,  this was stated for the WLP, but the same proof also gives it for the SLP.
\end{remark}

For the remainder of this section we will {\em  assume that $k$ has characteristic 0, }
unless explicitly mentioned otherwise.   We have seen that in one or
two variables, we always have the WLP (and even SLP).  Turning to
the next case, the following seemingly simple result in fact has a
very intricate and long proof.  It illustrates the subtlety of these
problems.

\begin{theorem}[\cite{BMMNZ}, Theorem 6.2]
  \label{thm:WLP-type-2}
A level artinian monomial algebra of type 2 in three variables has the WLP.
\end{theorem}

Of course this has the following consequence.

\begin{corollary}
A pure $O$-sequence of type 2 in three variables is unimodal.
\end{corollary}

The monograph \cite{BMMNZ} gave a careful study of families of level
artinian monomial algebras that fail the WLP.  As a consequence, we
have the following conclusion.

\begin{theorem}[\cite{BMMNZ}]
If $R = k[x_1,\dots,x_r]$ and $R/I$ is a level artinian monomial
algebra of type $t$, then, for all $r$ and $t$, examples exist where
the WLP fails, except if:

\begin{itemize}
\item $r=1$ or 2;

\item $t=1$ (this is Theorem \ref{swrrr thm});

\item $r=3$, $t=2$ (this is Theorem \ref{thm:WLP-type-2}).
\end{itemize}

\end{theorem}

In particular, the first case where WLP can fail is when $r=3$ and
$t=3$.  This occurs, for instance, if $I = \langle x^3, y^3, z^3,
xyz \rangle$  (see \cite{BK-WLP}, Example 3.1).  Nevertheless, Boyle
has shown that despite the failure of the WLP, all level artinian
monomial algebras with $r=3$ and $t=3$ have {\em strictly unimodal}
Hilbert function (that is, in addition to being unimodal, once the
function decreases then it is strictly decreasing from that point
until it reaches zero):

\begin{theorem}[\cite{boyle}]
Any pure $O$-sequence  of {\em  type 3} in three variables is
strictly unimodal.
\end{theorem}

In more variables, the first case where the WLP can fail is when
$r=4$ and $t=2$.  Here again, Boyle has shown that nevertheless such
algebras have strictly unimodal Hilbert function:

\begin{theorem}[\cite{boyle2}]
Any pure $O$-sequence of type 2 in four variables is strictly
unimodal.
\end{theorem}

\noindent Since the WLP is not available in these cases, Boyle's
method is a classification theorem followed by a decomposition of
the ideals and a careful analysis of sums of Hilbert functions of
complete intersections.

However, there is no hope of such a result for all pure
$O$-sequences,  even when $r=3$:

\begin{theorem}[\cite{BMMNZ}]
Let $M$ be any positive integer and fix an integer $r \geq 3$.
Then there exists a pure $O$-sequence in $r$ variables which is
non-unimodal, having exactly $M$ maxima.
\end{theorem}

In view of the last two results, we have the following natural
question.

\begin{question}
What is  the smallest socle degree and (separately) the smallest
socle type $t$ for which non-unimodal pure $O$-sequences exist.
This is especially of interest when $r=3$.
\end{question}

In \cite{BZ}, Boij and Zanello produced a non-unimodal example with
$r=3$ and socle degree 12.   In \cite{BMMNZ}, for $r=3$, we produced
a non-unimodal example  for socle type $t = 14$. It was also shown
that pure $O$-sequences can fail unimodality iff the socle degree is
at least 4 (but one may need many variables for small socle degree).

It is also  natural to ask how things change when you remove
``monomial" and ask about artinian level algebras.  Some work in
progress by Boij, Migliore, Mir\'o-Roig, Nagel and Zanello indicates
that the behavior of such algebras from the point of view of the
Hilbert function can become  surprisingly worse, in the sense that
dramatic non-unimodality is possible even in early degrees, which
would violate Hausel's theorem (Theorem \ref{hausel}) for instance,
in the monomial case.


\section{Powers of linear forms}

In this section we always assume that $k$ has characteristic zero.
Note that $x_i$ is a linear form, and that if $L_1,\dots,L_n$ $(n
\geq r)$ are general linear forms, then without loss of generality
(by a change of variables) we can assume that $L_1 = x_1,\dots,L_r =
x_r$. Thus Theorem \ref{swrrr thm} is also a result about ideals
generated by powers of linear forms.  It says that in
$k[x_1,\dots,x_r]$, an ideal generated by  powers of $r$ general
linear forms has the WLP and the SLP.  It also leads to an
interesting connection to Fr\"oberg's conjecture, which we discuss
in Section \ref{froeberg section}.

\begin{question}
Which ideals generated by powers of general linear forms define
algebras that fail the WLP or SLP?
\end{question}

We saw in Theorem \ref{hmnw 2 var} that all such ideals (and  in
fact all homogeneous ideals) in two variables satisfy both the WLP
and the SLP.  More surprisingly, Schenck and Seceleanu showed a
similar result in three variables:

\begin{theorem}[\cite{ss}] \label{ss 3 vars}
Let $R = k[x,y,z]$, where char$(k) = 0$.     Let $I = \langle
L_1^{a_1},\dots,L_m^{a_m} \rangle$ be {\em any} ideal generated by
powers of linear forms.  Then $R/I$ has the WLP.
\end{theorem}

\vspace{.2in}

A shorter proof of this result is given in \cite{MMN-powers}.  One
reason that it is surprising is that the same is {\em not} true for
SLP.  For instance, if $I = \langle L_1^3,L_2^3,L_3^3,L_4^3 \rangle$
(where $L_i$ general in $k[x,y,z]$), then $(\times \ell^3)$ fails to
have maximal rank. The case of three variables acts as a bridge
case: we will see that for four or more variables, even WLP fails
very commonly. Some recent work in this area was motivated by the
following example  of Migliore, Mir\'o-Roig and Nagel:

\begin{example}[\cite{MMN}] \label{MMN orig ex} Let $r=4$.    Consider the ideal $I = \langle x_1^N,x_2^N,x_3^N,x_4^N,L^N \rangle$ for a general linear form $L$.
By computation using CoCoA,  $R/I$ fails the WLP, for
$N=3,\dots,12$.
\end{example}

There are some natural
questions arising from this example:

\begin{problem}

\begin{itemize}
\item Prove the failure of the WLP in Example \ref{MMN orig ex} for
 all $N \geq 3$.

\medskip

\item What happens for mixed powers?

\medskip

\item What happens for almost complete intersections, that is, for $r+1$ powers of general linear forms in $r$ variables when $r>4$?

\medskip

\item What about more than $r+1$ powers of  general linear forms?

\end{itemize}
\end{problem}

This example motivated two different projects at the same time: by Migliore, Mir\'o-Roig, Nagel \cite{MMN-powers}  and by Harbourne, Schenck, Seceleanu \cite{HSS}.   Both of these papers used the dictionary between ideals of powers of general linear forms and ideals of fat points in projective space, provided by the following important result of Emsalem and Iarrobino:

\begin{theorem}[\cite{EI}]
Let
\[
\langle L_1^{a_1},\dots,L_n^{a_n} \rangle \subset k[x_1,\dots,x_r]
\]
 be an ideal generated by powers of $n$  linear forms.  Let $\wp_1,\dots,\wp_n$ be the ideals of $n$ the points in $\mathbb P^{r-1}$ corresponding to the linear forms.  Then for any integer $j \geq \max \{ a_i\}$,
\[
\dim_k [R/\langle L_1^{a_1},\dots,L_n^{a_n} \rangle ]_j = \dim_k \left [ \wp_1^{j-a_1+1} \cap \cdots \cap \wp_n^{j-a_n+1} \right ]_j.
\]

\end{theorem}

One important difference between the two papers is that \cite{HSS} assumed that the powers are uniform, and usually that the powers are ``large enough.''  Usually they allow more than $r+1$ forms.  On the other hand, most of the results in \cite{MMN-powers} allow mixed powers.  We quote some of the results of these two papers.  Together they form a nice start to an interesting topic. The conjectures listed later indicate that more work is to be done!

\begin{theorem} [\cite{HSS}] \label{hss th1}
 Let
 \[
 I = \langle L_1^t,\dots,L_n^t \rangle \subset k[x_1,x_2,x_3,x_4]
 \]
  with $L_i \in R_1$ generic.    If $n \in \{ 5,6,7,8 \}$,  then the WLP fails, respectively, for $t \geq \{ 3, 27, 140, 704 \}$.
 \end{theorem}

\begin{theorem}[\cite{HSS}] \label{hss th2}
For
\[
I = \langle L_1^t,\dots,L_{2k+1}^t \rangle \subset R = k[x_1,\dots,x_{2k}]
\]
 with $L_i$ generic linear forms, $k\geq 2$ and $t \gg 0$, $R/I$ fails the WLP.
\end{theorem}

\noindent (See also Theorem \ref{even uniform}  below.)  The following result gives the most complete picture to date, about the case of four variables, when the exponents are not assumed to be uniform and the ideal is assumed to be an almost complete intersection (i.e. the number of minimal generators is one more than the number of variables).  It summarizes several theorems in \cite{MMN-powers}, Section 3, and we refer to that paper for the more detailed individual statements.

\begin{theorem} [Four variables, \cite{MMN-powers}]  \label{MMN 4 vars}
 Let
 \[
 I = \langle L_1^{a_1},L_2^{a_2},L_3^{a_3},L_4^{a_4},L_5^{a_5} \rangle \subset R = k[x_1,x_2,x_3,x_4],
 \]
  where all $L_i$ are generic.
Without loss of generality assume that  $a_1 \leq a_2 \leq a_3 \leq a_4 \leq a_5$.  Set
\[
\lambda =
\left \{
\begin{array}{ll}
\displaystyle \frac{a_1+a_2+a_3+a_4}{2} - 2 & \hbox{if $a_1+a_2+a_3+a_4$ is even} \\ \\
\displaystyle \frac{a_1+a_2+a_3+a_4-7}{2} & \hbox{if $a_1+a_2+a_3+a_4$ is odd.}
\end{array}
\right.
\]

\begin{itemize}

\item[(a)] If $a_5 \geq \lambda$
 then $R/I$ has the WLP.

\item[(b)] If $a_1 = 2$ then $R/I$ has the WLP.

\item[(c)] Most other cases (explicitly described in terms of  $a_1,a_2,a_3,a_4$) are proven to fail the WLP.

\item[(d)] For the few open cases, experimentally sometimes the WLP holds and sometimes not.

\end{itemize}
\end{theorem}

\noindent Notice that the case where all the $a_i$ are equal and at least $3$ is contained in Theorem \ref{hss th1}.

In more than four variables, it becomes progressively more difficult to obtain results for mixed powers.  We have the following partial result.

\begin{theorem} [Five variables, almost uniform powers, \cite{MMN-powers}] \label{five var}
 Assume $r=5$.   Let $L_1,\dots,L_6$ be  general linear forms.    Let $e \geq 0$ and
\[
I = \langle L_1^d,L_2^d,L_3^d,L_4^d,L_5^d,L_6^{d+e} \rangle.
\]

\begin{itemize}

\item[(a)]  If $e=0$ then $R/I$ fails the WLP if and only if $d > 3$.

\medskip

\item[(b)] If $e \geq 1$ and $d$ is odd then $R/I$ has the WLP if and only if $e \geq \frac{3d-5}{2}$.

\medskip

\item[(c)] If $e \geq 1$ and $d$ is even then $R/I$ has the WLP if and only if $e \geq \frac{3d-8}{2}$.
\end{itemize}
\end{theorem}

We also have the following improvement of Theorem \ref{hss th2}, which has the additional assumption that $t \gg 0$.

\begin{theorem} [Even number of variables, uniform powers, \cite{MMN-powers}] \label{even uniform}

Let
\[
I = \langle L_1^t,\dots,L_{2k+1}^t \rangle \subset R = k[x_1,\dots,x_{2k}]
\]
 with $L_i$ generic linear forms and $k\geq 2$.
  Then $R/I$ fails the WLP if and only if $t>1$.

\end{theorem}

\noindent (The case $k=2$ is contained in Theorem \ref{MMN 4 vars}.)

What about an odd number of variables?    Here is a result for seven variables:

\begin{theorem} [\cite{MMN-powers}] \label{7vars}
 Let
 \[
 I = \langle L_1^t,\dots,L_8^t \rangle \subset k[x_1,\dots,x_7],
 \] where $L_1,\dots,L_8$ are general linear forms.

 \medskip

\begin{itemize}

\item If $t=2$ then $R/I$ has the WLP.

\medskip

\item  If $t \geq 4$ then $R/I$ fails the WLP.

\end{itemize}
\end{theorem}

Interestingly, for $t=3$, CoCoA \cite{cocoa} says that the WLP fails, but we do not have a proof.  We can believe a computer that says that the WLP {\em holds}, but otherwise we have to be skeptical about whether its choice of forms was ``general enough."

For  these results, sometimes it was necessary to prove {\em failure of surjectivity} (when $h_{i-1} \geq h_i$ in the relevant degrees), sometimes {\em failure of injectivity} (when $h_{i-1} \leq h_i$), and sometimes we had to show that {\em the WLP does hold}.    These present quite different challenges.

After making the translation to fat points, as described above, the first difficulty is to determine the degrees where WLP fails.  Then, it was necessary to find the dimension of a linear system of surfaces in a suitable projective space vanishing to prescribed multiplicity at a general set of suitably many points.  To do this, in \cite{MMN-powers} Cremona transformations and work of Dumnicki 2009, Laface-Ugaglia 2006, and De Volder-Laface 2007 were used as main tools, plus ad hoc methods.  These Cremona transformation results are central to the proofs in \cite{MMN-powers}.

Much remains to be shown on this topic.  Here are two conjectures from \cite{HSS} and \cite{MMN-powers}.

\begin{conjecture} [\cite{HSS}]
For $I = \langle L_1^t,\dots,L_n^t \rangle \subset R = [x_1,\dots,x_r]$ with $L_i \in R_1$ generic and $n \geq r+1 \geq 5$, the WLP fails for all $t \gg 0$.
\end{conjecture}

\begin{conjecture}[\cite{MMN-powers}]

Let $R = k[x_1,\dots,x_{2n+1}]$.  Let $L_1,\dots,L_{2n+2}$ be  general linear forms and $I = \langle L_1^d,\dots,L_{2n+1}^d,L_{2n+2}^d \rangle$.

\bigskip

\begin{itemize}

\item If $n =3$ and $d=3$ then $R/I$ fails the WLP. (This is the only open case in Theorem \ref{7vars}.)

\bigskip

\item If $n \geq 4$ then $R/I$ fails the WLP if and only if $d>1$.

\end{itemize}

\end{conjecture}

These conjectures are supported by a great deal of computer evidence, using \cocoa \cite{cocoa} and {\tt Macaulay2} \cite{macaulay}.


\section{Connection between Fr\"oberg's conjecture and the WLP} \label{froeberg section}

In this section we continue to assume that our field has characteristic zero.  Closely related to the SLP is the so-called {\em maximal rank property (MRP)}, which just replaces $\ell^d$ by a general form of degree $d$ in Definition \ref{def of wlp slp}.  Nevertheless, it is known that  the MRP does not imply the SLP.  See \cite{MM} and \cite{ZZ} for some connections between these two properties.

One way of stating Fr\"oberg's conjecture is as follows.
\begin{conjecture}[Fr\"oberg]
Any ideal of general forms has the MRP.
More precisely, fix positive integers $a_1,\dots,a_s$ for some $s>1$.  Let $F_1,\dots,F_s \subset R = k[x_1,\dots,x_r]$ be general forms of degrees $a_1,\dots,a_s$ respectively and let $I = \langle F_1,\dots,F_s \rangle$.  Then for each $i$, $2 \leq i \leq s$, and for all $t$, the multiplication by $F_i$ on $R/\langle F_1,\dots,F_{i-1} \rangle$ has maximal rank, from degree $t-a_i$ to degree $t$.  As a result, the Hilbert function of $R/I$ can be computed inductively.
\end{conjecture}

This conjecture is known to be true in two variables. This follows, for example, from Theorem \ref{hmnw 2 var}.  In three variables it was shown to be true by Anick \cite{anick}.
In this section, we explore the following natural questions.

\begin{question}
What is the Hilbert  function of an ideal generated by powers of general linear forms of degrees $a_1,\dots,a_n$?    In particular, is it the same as the Hilbert function predicted by Fr\"oberg?  What, if any, is the connection to the WLP?
\end{question}

\noindent  Theorem \ref{swrrr thm} says that when $n = r+1$, the answer to the second question is yes.

The fact that the answer is often ``no" for $n = r+2$ was first observed by Iarrobino \cite{iarrobino}.  Chandler \cite{chandler-JA}, \cite{chandler-Con} also gave some results in this direction.  Concerning the connection to the WLP, the following result of Migliore, Mir\'o-Roig and Nagel gives a partial answer.
 
\begin{proposition}[\cite{MMN-powers}]
  \label{prop:Lefschetz-prop}
\begin{itemize}
  \item[(a)]  If Fr\"oberg's conjecture is true for all ideals generated by general forms in $r$ variables, then all ideals generated by general forms in $r+1$ variables have the WLP.

  \item[(b)]    Let $R = k[x_1,\dots,x_{r+1}]$, let $\ell \in R$ be a general linear form, and let $S =R/\langle \ell \rangle \cong  k[x_1,\dots,x_{r}]$.  Fix positive integers  $s, d_1,\dots,d_s, d_{s+1}$.  Let $L_1,L_2,\dots,L_s,L_{s+1} \in R$ be  linear forms.  Denote by $\bar {\hbox{\ \ }}$ the restriction from $R$ to $S \cong R/\langle \ell \rangle$.  Make the following assumptions:

\begin{itemize}
\item[(i)] The ideal $I = \langle L_1^{d_1},\dots,L_s^{d_s} \rangle$ has the WLP.

  \item[(ii)] The multiplication $\times \bar L_{s+1}^{d_{s+1}} : [S/\bar I]_{j-d_{s+1}} \rightarrow [S/\bar I]_j$ has maximal rank.
     \end{itemize}

 Then  $R/\langle L_1^{d_1},\dots,L_{s+1}^{d_{s+1}} \rangle$ has the WLP.

\end{itemize}
\end{proposition}

\begin{remark}
\begin{itemize}
\item[(a)] Part of this result was in fact contained in the paper \cite{MM} of Migliore and Mir\'o-Roig.  It was used there  to show that any ideal of general forms in $k[x_1,x_2,x_3,x_4]$ satisfies the WLP, because Anick \cite{anick} had shown much earlier that any ideal of general forms in $k[x_1,x_2,x_3]$ satisfies Fr\"oberg's conjecture.

\medskip

\item[(b)] It was shown in \cite{MMN-powers} that this result also leads to a short proof of Theorem \ref{ss 3 vars}.  The point is that the restriction of such ideals correspond to an ideal in $k[x,y]$, and in characteristic zero all such ideals have the SLP by Theorem \ref{hmnw 2 var}.
\end{itemize}
\end{remark}

\medskip

The following corollary was also shown in \cite{MMN-powers}:

\begin{corollary}[\cite{MMN-powers}] \label{appl}
  Assume the characteristic is zero.  Let $R = k[x_1,\dots,x_{r+1}]$, let $\ell \in R$ be a general linear form, and let $S =R/\langle \ell \rangle \cong  k[x_1,\dots,x_{r}]$.  For integers $d_1,\dots,d_{r+2}$, if an ideal of the form $\langle L_1^{d_1},\dots,L_{r+2}^{d_{r+2}} \rangle \subset R$ of powers of general linear forms fails to have the WLP then an ideal of powers of general linear forms $\langle \bar L_1^{d_1},\dots,\bar L_{r+2}^{d_{r+2}} \rangle \subset S$ fails to have the Hilbert function predicted by Fr\"oberg's conjecture.
 \end{corollary}

Thus the results in the previous section give additional insight to the observation of Iarrobino \cite{iarrobino} and Chandler \cite{chandler-JA}, \cite{chandler-Con} that when $n=r+2$, there are many cases when an ideal of powers of general linear forms does not have the same Hilbert function as that predicted by Fr\"oberg for general forms.  Since Theorem \ref{MMN 4 vars} covers almost all possible choices of exponents, it gives a much more complete answer to the question of exactly which powers of five general linear forms in three variables fail to have the Fr\"oberg-predicted Hilbert function, contrasting with the result of Anick which says that an ideal of general forms of any fixed degrees in three variables does have the predicted Hilbert function.  Theorems \ref{five var} and \ref{even uniform} provide new partial answers (via Corollary \ref{appl}) in the case of more variables.

\begin{example}
Let $R = k[x_1,x_2,x_3,x_4]$.  Let $L_1, L_2, L_3, L_4, L_5$ and $\ell$ be  general linear forms.  Let $S = R/\langle L \rangle \cong k[x,y,z]$.  Let $I = \langle L_1^3,L_2^3,L_3^3, L_4^3, L_5^3 \rangle$.  (The smallest case in Example \ref{MMN orig ex}  above.)
The Hilbert function of $R/I$ is $(1,4,10, 15, 15, 6)$.
We have
\[
[R/I]_3 \stackrel{\times \ell}{\longrightarrow} [R/I]_{4} \rightarrow [R/(I,\ell)]_{4} \rightarrow 0.
\]
We saw that WLP fails, and in fact
\[
\dim [R/(I,\ell)]_4 = 1.
\]
Notice that $R/(I,\ell) \cong S/J$, where $J$ is the ideal of cubes of five general linear forms in $k[x,y,z]$.
Thus $\dim [S/J]_4 = 1$.

On the other hand, let $K$ be the ideal of five general cubics in $S$.  Fr\"oberg predicts (and Anick proves) that $\dim[S/K]_4 = 0$.  Thus $J$ does not have the Hilbert function predicted by Fr\"oberg.

In fact, whenever we prove that an ideal of $n$ powers of general linear forms fails the WLP (for specified exponents), then for some subset of these powers of general linear forms, the same number and powers of general linear forms in one fewer variable fails to have Fr\"oberg's predicted Hilbert function.
\end{example}

\section{Positive characteristic and enumerations} \label{pos char}

Considering Theorem \ref{swrrr thm} again, we saw in Lemma
\ref{lem:ci-pos-char} that the assumption on the characteristic
of the base field cannot be omitted.

\begin{question}
   \label{q:ci-pos-char}
What happens in Theorem \ref{swrrr thm} if we allow the
characteristic to be positive?
\end{question}

Actually, investigating the dependence of the WLP on the
characteristic makes sense whenever the ideal can be defined
over the integers. This applies to all monomial ideals. In fact,
in this case one has the following result.

\begin{proposition}[\cite{CN2}, Lemma 2.6]
  \label{prop:char-zero-to-pos} Let $I \subset R$ be a monomial
  ideal. If $R/I$ has the WLP when $\chara (k) = 0$, then $R/I$
  has the WLP whenever $\chara (k)$ is sufficiently large.
\end{proposition}

The proof is based on two observations that have their origin in
\cite{MMN}. For a monomial ideal, one can check the WLP by
considering the specific linear form $\ell = x_1 + \cdots +
x_r$. Thus, the maximal rank property of the multiplications by
$\ell$ is governed by integer matrices. Their determinants have
only finitely many prime divisors if they do not vanish.

It also follows that $R/I$ fails the WLP in {\em every} positive
characteristic if it fails the WLP in characteristic zero.

Proposition \ref{prop:char-zero-to-pos} motivates the following
problem.

\begin{question}
  \label{quest:char-WLP-fails} Let $I \subset R$ be a monomial
ideal such that $R/I$ has the WLP when $\chara k = 0$. What are
the (finitely many) field characteristics such that $R/I$ fails
the WLP?
\end{question}

This turns out to be a rather subtle problem. It was first
considered in \cite{MMN} in the case of certain almost complete
intersection in three variables. Recall that a monomial almost
complete intersection in three variables is an ideal of the form
\begin{equation}
  \label{eq-aci}
I =   I_{a,b,c,\alpha,\beta,\gamma} =
    \langle x^a, y^b, z^c, x^\alpha y^\beta z^\gamma \rangle.
\end{equation}
If the syzygy bundle of $I$ is not semistable or its first Chern
class is not divisible by three, then $R/I$ has the WLP in
characteristic zero (see \cite{BK-WLP} and \cite{MMN}).  However, if
the syzygy bundle satisfies both conditions, then deciding the WLP
is more difficult and very subtle on the one side. On the other
side, the investigations in this case have brought to light
surprising connections to combinatorial problems.

In fact, if the syzygy bundle  of $I$ is semistable and its first Chern class is divisible by three, then $R/I$ has the WLP if and only
if the multiplication by $\ell$ in a certain degree is an
isomorphism or, equivalently, a certain integer square matrix has a
non-vanishing determinant. This has been first observed in the
special case, where $R/I$ is level, in \cite{MMN} and then for
arbitrary almost complete intersections in \cite{CN2}. The first
connection to combinatorics was made by Cook  and the second author
in Section 4 of \cite{CN}. There it was observed  that the
determinant  deciding the WLP for certain families of monomial
almost complete intersections is the number of lozenge tilings of
some hexagon, which is given by a formula of MacMahon.  Lozenge
tilings of a hexagon are in bijection to other well-studied
combinatorial objects such as, for example, plane partitions and
families of non-intersecting lattice paths.

Independently of \cite{CN} but subsequent to it, Li and Zanello
studied the WLP in the case of the complete intersections $R/\langle
x^a,y^b,z^c \rangle$ in \cite{LZ}, and they also related MacMahon's
numbers of plane partitions to the failure of the WLP:

\begin{theorem}[\cite{LZ}]
For any given positive integers $a,b,c$, the number of plane
partitions contained inside an $a \times b \times c$ box is
divisible by a  prime $p$ if and only if the algebra $k[x,y,z]/
\langle x^{a+b},y^{a+c},z^{b+c} \rangle$ fails to have the WLP when
$\hbox{char}(k) = p$.
\end{theorem}

\noindent (This connection is already implicitly contained in
\cite{CN}, although it was only made explicit in the proof of
Corollary 6.5 in \cite{CN2}.) Next, Chen, Guo, Jin, and Liu
\cite{CGJL}, explained \emph{bijectively}  the result by Li and
Zanello for complete intersections.  Both \cite{LZ} and \cite{CN}
have been substantially extended in \cite{CN2}.  Here the bijective
approach of \cite{CGJL} was extended to almost complete
intersections, and further relations between the presence of the WLP
and difficult counting problems in combinatorics have been given.
In the remainder of this section, we give an overview of some of the
results of \cite{CN2} which illustrate this fascinating connection.

We focus on the most difficult case, in which the presence of the
WLP is a priori not even known  in characteristic zero, that is, we
assume that the syzygy bundle of the almost complete intersection $I
= I_{a,b,c,\alpha,\beta,\gamma} = \langle x^a, y^b, z^c, x^\alpha
y^\beta z^\gamma \rangle$ is semistable in characteristic zero and
its first Chern class is divisible by three.  By \cite[Proposition
3.3]{CN2}, this is exactly true if and only if the following
conditions  are all satisfied:
\begin{itemize}
        \item[(i)] $s := \frac{1}{3}(a + b + c + \alpha + \beta + \gamma)-2$ is an integer,
        \item[(ii)] $0 \leq M$,
        \item[(iii)] $0 \leq A \leq \beta  + \gamma$,
        \item[(iv)] $0 \leq B \leq \alpha + \gamma$, and
        \item[(v)] $0 \leq C \leq \alpha + \beta$,
\end{itemize}
where
\begin{equation*}
        \begin{split}
            A :=& s+2-a, \\
            B :=& s+2-b, \\
            C :=& s+2-c, \mbox{ and}\\
            M :=& s + 2 - (\alpha + \beta + \gamma).
        \end{split}
    \end{equation*}

\begin{figure}[!ht]
        \includegraphics[scale=0.667]{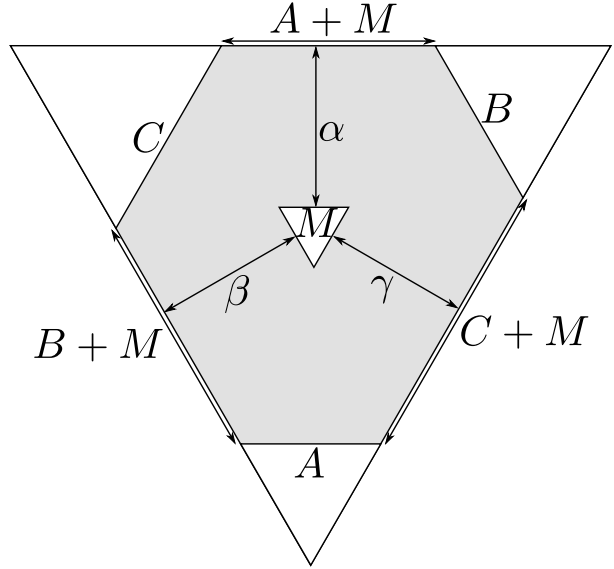}
        \caption{Punctured hexagon $H_{a,b,c,\alpha,\beta,\gamma}$ (shadowed) associated to $I_{a,b,c,\alpha,\beta,\gamma}$}
        \label{fig:punctured-hexagon}
\end{figure}

The above conditions have a geometric meaning. In fact, due  to
Theorem 4.1 in \cite{CN2},  they  guarantee that $I$ can be related
to a hexagonal region with a hole, which is called the {\em
punctured hexagon} $H = H_{a,b,c,\alpha,\beta,\gamma}$ associated to
$I = I_{a,b,c,\alpha,\beta,\gamma}$ (see
Figure~\ref{fig:punctured-hexagon}).

There are two square matrices that govern the WLP of the ideal  $I$.
In fact, $I$ has the WLP if and only if the multiplication $[R/I]_s
\stackrel{\times \ell}{\longrightarrow} [R/I]_{s+1}$  is bijective
or, equivalently, $[R/(I, \ell)]_{s+1} = 0$. The latter condition
means that a certain $(C+M) \times (C+M)$ matrix, $N =
N_{a,b,c,\alpha,\beta,\gamma}$, with binomial coefficients as
entries is regular. The above multiplication map can be described by
a much larger zero-one square matrix,  $Z =
Z_{a,b,c,\alpha,\beta,\gamma}$. The mentioned equivalence implies
that the determinants of $N$ and $Z$ have the same prime divisors.
However, much more is true. Both determinants have the same absolute
value, which has combinatorial interpretations.

\begin{theorem}[\cite{CN2}, Theorems 4.4, 4.6, and 5.4]
  \label{thm:aci-enumerations}
Adopt the above assumptions. Then the following
conditions are equivalent:
\begin{itemize}
  \item[(a)] $I_{a,b,c,\alpha,\beta,\gamma}$ has the WLP if the
  characteristic of the base field $k$ is $p \geq 0$.

  \item[(b)] $p$ does not divide the enumeration
  $|\det N_{a,b,c,\alpha,\beta,\gamma}|$ of signed lozenge
  tilings of the associated punctured hexagon
  $H_{a,b,c,\alpha,\beta,\gamma}$.

  \item[(c)] $p$ does not divide the enumeration
  $|\det Z_{a,b,c,\alpha,\beta,\gamma}|$ of signed perfect
  matchings of the bipartite graph associated to
  $H_{a,b,c,\alpha,\beta,\gamma}$.
\end{itemize}
In particular, $|\det N_{a,b,c,\alpha,\beta,\gamma}| = |\det
Z_{a,b,c,\alpha,\beta,\gamma}|$.
\end{theorem}

\begin{figure}[!ht]
    \begin{minipage}[b]{0.49\linewidth}
        \centering
        \includegraphics[scale=0.45]{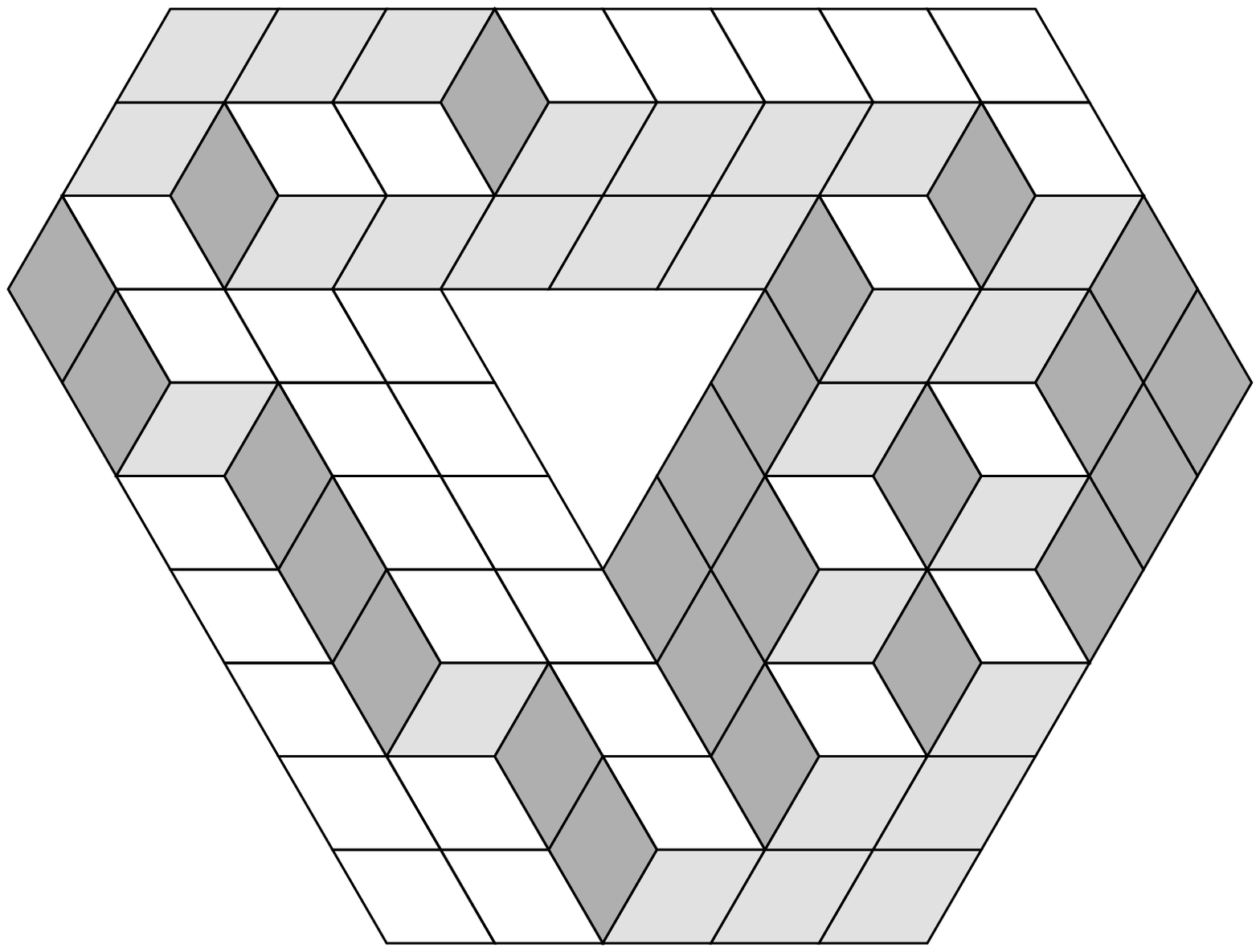}\\
        {\em Hexagon tiling by lozenges}
    \end{minipage}
    \begin{minipage}[b]{0.49\linewidth}
        \centering
        \includegraphics[scale=0.45]{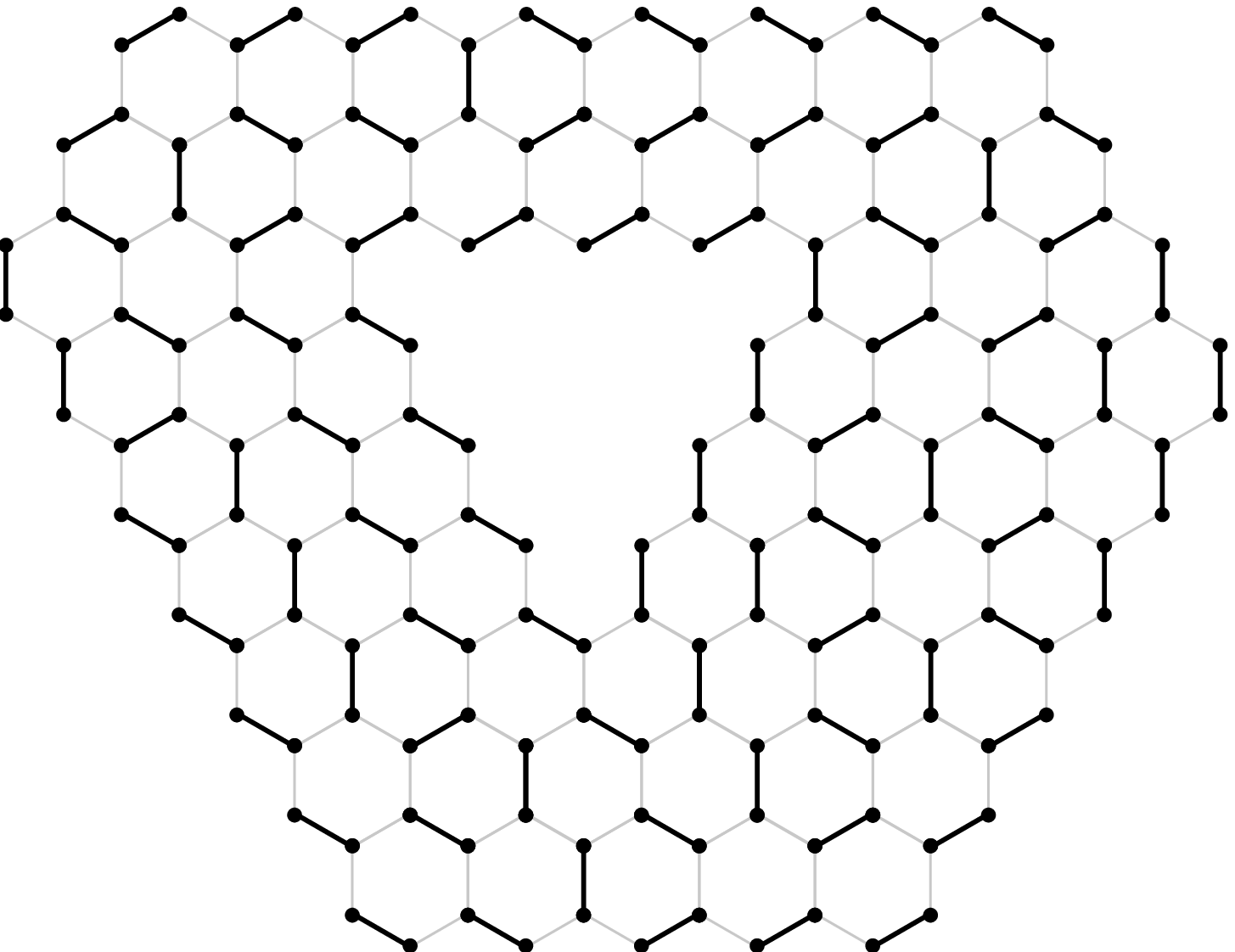}\\
        {\em Perfect matching of edges}
    \end{minipage}
    \caption{A lozenge tiling and its associated perfect matching}
    \label{fig:hex-bip}
\end{figure}

A lozenge is a rhombus with unit side-lengths and angles of
$60^\circ$ and $120^\circ$.  Lozenges have also been called
calissons and diamonds in the literature. A perfect matching of a
graph is a set of pairwise non-adjacent edges such that each vertex
of the graph is matched. We refer to \cite{CN2} for more details, in
particular for assigning the signs,  though Figure \ref{fig:hex-bip}
below  indicates an associated lozenge tiling and a perfect
matching.

Theorem \ref{thm:aci-enumerations} has been used to establish the
WLP of $I_{a,b,c,\alpha,\beta,\gamma}$ in many new cases. The
results also lend further evidence to a conjectured characterization
of the presence of the WLP of $I_{a,b,c,\alpha,\beta,\gamma}$ in
case $I_{a,b,c,\alpha,\beta,\gamma}$ is level that has been proposed
in \cite{MMN}.

The determinants occurring in Theorem \ref{thm:aci-enumerations} can
be rather big.

\begin{example}
  \label{ex:det}
Consider the ideal
\[
I = \langle x^{14}, y^{21}, z^{25}, x^2 y^9 z^{13}\rangle.
\]
Then the absolute value of the corresponding determinants is (see
\cite{CN}, Remark 4.8)
\[
2 \cdot 3^2 \cdot 5^3 \cdot 11^4 \cdot 13^5 \cdot 19 \cdot 23^3
\cdot 29 \cdot 5011.
\]
Hence $R/I$ fails the WLP if and only if the characteristic of $k$
is  any  of the nine  listed prime divisors.
\end{example}

In the situation of Theorem \ref{thm:aci-enumerations},  the
presence of the WLP in characteristic zero can also be read off from
the splitting type of the syzygy bundle. In fact,
$I_{a,b,c,\alpha,\beta,\gamma}$ has the WLP if and only if its
syzygy bundle has splitting type $(s+2, s+2, s+2)$ (see \cite{CN2},
Theorem 9.9).

In \cite{CN2}, explicit formulae for the enumerations appearing in
Theorem \ref{thm:aci-enumerations} are derived in various cases.
However, even then determining the prime divisors of the
enumerations can be challenging. In fact, this problem is even open
in the special case of monomial complete intersections. Though,
recently, there has been progress in the case, where the generators
all have the same degree. Brenner and Kaid \cite{BK2-char p} gave an
explicit description of when $R/\langle x^d,y^d,z^d \rangle$  has
the WLP in terms of $d$ and the characteristic $p$.  In particular, they proved a conjecture of \cite{LZ} for the case $p=2$. This latter result is stated very concisely:

\begin{theorem}\cite{BK2-char p}
The algebra $k[x, y, z]/\langle x^d, y^d, z^d \rangle$ has the WLP
in $\hbox{char}(k) = 2$ if and only if $d = \lfloor \frac{2^n +1}{3}
\rfloor$ for some positive integer $n$.
\end{theorem}

\noindent The approach of \cite{BK2-char p} was via a theorem of Han
computing the syzygy gap for an ideal of the form $\langle
x^d,y^d,(x+y)^d \rangle$ in $k[x,y]$. The analogous result in the
case of more variables, that is, for $I = \langle x_1^d,\ldots,x_n^d
\rangle$ ($n \ge 4$), has been obtained by Kustin and Vraciu in
\cite{KV}. Independently, Cook made progress in deciding the
Lefschetz properties of more general monomial complete intersections
in positive characteristic (see \cite{C}), addressing Question
\ref{q:ci-pos-char} (see also \cite{lindsey}, Lemma 5.2, for a
result in two variables).

In a different direction, Kustin, Rahmati and Vraciu \cite{KRV}
showed that $A = R/\langle x^d,y^d,z^d \rangle$ has the WLP in
characteristic $p \neq 2$ if and only if its residue field has finite
projective dimension as an $A$-module.




\begin{thebibliography}{999}

\bibitem{amasaki}

M.\ Amasaki, {\em The weak Lefschetz property for Artinian graded
rings and basic sequences}, Preprint, 2011.

\bibitem{anick} D.\ Anick, {\em Thin algebras of embedding dimension three}, J.\ Algebra {\bf 100} (1986), 235--259.

\bibitem{BI}  D.\ Bernstein and A.\ Iarrobino: {\it A
nonunimodal graded Gorenstein Artin algebra in codimension five},
Comm.\ Algebra {\bf 20} (1992), 2323--2336.

\bibitem{Bo} M.\ Boij: {\it Graded Gorenstein Artin algebras
whose Hilbert functions have a large number of valleys}, Comm.\
Algebra {\bf 23} (1995),  97--103.

\bibitem{BL} M.\ Boij and D.\ Laksov: {\em Nonunimodality of
graded Gorenstein Artin algebras}, Proc.\ Amer.\ Math.\ Soc.\ {\bf 120}
(1994), 1083-1092.

\bibitem{BZ} M.\ Boij and F. Zanello, {\em Level algebras with bad properties}, Proc. Amer. Math. Soc. {\bf  135} (2007), no. 9, 2713--2722.

\bibitem{BMMNZ} M.\ Boij, J.\ Migliore, R.\ M.\ Mir\'o-Roig,
U.\ Nagel and F.\ Zanello, ``The shape of a pure $O$-sequence,''
Mem.\ Amer.\ Math.\  Soc.\ (to appear).

\bibitem{boyle} B.\ Boyle, {\em The unimodality of pure $O$-sequences of type three in three variables}, preprint  2011.

\bibitem{boyle2} B.\ Boyle, {\em The unimodality of pure $O$-sequences of type two in four variables}, preprint  2011.

\bibitem{BK-WLP} H.\ Brenner and A.\ Kaid, {\em Syzygy bundles on $\mathbb P^2$ and the Weak Lefschetz Property}, Illinois J.\ Math.\ {\bf 51} (2007), 1299--1308.

\bibitem{BK2-char p} H.\ Brenner and A.\ Kaid, {\em A note on the weak Lefschetz property of monomial
complete intersections in positive characteristic}, Collect.\ Math.\ {\bf 62} (2011), 85--93.

\bibitem{chandler-JA} K.\ Chandler, {\em The geometric interpretation of Fr\"oberg-Iarrobino conjectures on infinitesimal neighborhoods of points in projective space}, J.\ Algebra {\bf 286} (2005), 421--455.

\bibitem{chandler-Con} K.\ Chandler, {\em Examples and counterexamples on the conjectured Hilbert function of multiple points}, in: ``Algebra, Geometry and Their Interactions, Contemp.\ Math., 448, Amer.\ Math.\ Soc., Providence, RI, 2007, 13--31.

\bibitem{CGJL} C.\ Chen, A.\ Guo, X.\ Jin and G.\ Liu,
{\em Trivariate monomial complete intersections and plane
partitions}, J.\ Commut.\ Algebra (to appear).

\bibitem{cocoa}   CoCoATeam, {\em
  CoCoA: a system for doing
     Computations in Commutative Algebra},
  Available at \url{http://cocoa.dima.unige.it}.

\bibitem{C} D.\ Cook II, {\em The Lefschetz properties of monomial
complete intersections in positive characteristic}, in preparation.

\bibitem{CN} D.\ Cook II and U.\ Nagel, {\em The Weak Lefschetz Property,
monomial ideals, and lozenges}, Illinois J.\ Math.\ (to appear).

\bibitem{CN-JPPA} D.\ Cook II and U.\ Nagel, {\em Hyperplane sections and the subtlety of the Lefschetz properties},  J.\ Pure Appl.\ Algebra {\bf 216} (2012),  108–-114.

\bibitem{CN2} D.\ Cook II and U.\ Nagel, {\em Enumerations deciding
the Weak Lefschetz Property}, Preprint, 2011.

\bibitem{EL} L.\ Ein and R.\  Lazarsfeld, {\em  Syzygies and Koszul
cohomology of smooth projective varieties of arbitrary dimension},
Invent.\ Math.\ {\bf 111} (1993), 51-–67.

\bibitem{EI} J.\ Emsalem and A.\ Iarrobino, {\em Inverse system of a symbolic power $I$}, J.\ Algebra {\bf 174} (1995), 1080-1090.

\bibitem{macaulay} D.\ Grayson and M.\ Stillman, {\em Macaulay2, a
software system for research in algebraic geometry}, Available
at  {\tt http://www.math.uiuc.edu/Macaulay2/}.

\bibitem{HSS} B.\ Harbourne, H.\ Schenck and A.\ Seceleanu,
{\em Inverse systems, Gelfand-Tsetlin patterns and the Weak
Lefschetz  Property}, J.\ London Math.\ Soc.\ (to appear).

\bibitem{HMMNWW} T. Harima, T. Maeno, H. Morita, Y. Numata, A. Wachi and J. Watanabe, ``The Lefschetz properties," manuscript available at {\tt http://www.stat.t.u-tokyo.ac.jp/~numata/tex/2010/hmmnww/} (2011).

\bibitem{HMNW} T.\ Harima, J.\ Migliore, U.\ Nagel, and J.\
Watanabe, {\em The Weak and Strong     Lefschetz Properties for
artinian $K$-Algebras}, J.\ Algebra {\bf 262} (2003), 99--126.

\bibitem{Ha} T.\ Hausel: {\em Quaternionic geometry of matroids}, Cent.\ Eur.\ J.\ Math.\ {\bf 3} (2005), no.\ 1, 26--38.\

\bibitem{HP} J.\ Herzog and D.\ Popescu, {\em The strong Lefschetz property and simple extensions}, preprint. Available on the arXiv at  http://front.math.ucdavis.edu/0506.5537.

\bibitem{iarrobino} A.\ Iarrobino, {\em Inverse system of a symbolic power III: thin algebras and fat points}, Compos.\ Math.\ {\bf 108} (1997), 319--356.

\bibitem{ikeda} H.\ Ikeda, {\em Results on Dilworth and Rees numbes of artinian
local rings}, Japan.\ J.\ Math.\ {\bf 22} (1996), 147--158.

\bibitem{KRV} A.\ Kustin, H.\ Rahmati and A.\ Vraciu, {\em The resolution of the bracket powers of the maximal ideal in a diagonal hypersurface ring}, Preprint,  2010.

\bibitem{KV} A.\ Kustin and A.\ Vraciu, {\em The Weak Lefschetz Property for monomial complete intersections in positive characteristic}, in preparation.

\bibitem{LZ} J.\ Li and F.\ Zanello, {\em Monomial complete intersections,
the Weak Lefschetz Property and plane partitions},  Discrete Math.
{\bf 310}  (2010),  no. 24, 3558–-3570.

\bibitem{lindsey}
M.\ Lindsey, {\em A class of Hilbert series and the strong Lefschetz
property}, Proc.\ Amer.\ Math.\ Soc.\ {\bf 139} (2011), 79–-92.

\bibitem{MM} J.\ Migliore and R.\ Mir\'o-Roig, {\em Ideals of general
forms and the ubiquity of the Weak Lefschetz property}, J.\ Pure
Appl.\ Algebra {\bf 182} (2003), 79--107.

\bibitem{MMN} J.\ Migliore, R.\ Mir\'o-Roig and U.\ Nagel, {\em
Monomial ideals,  almost complete intersections and the Weak
Lefschetz  Property}, Trans.\ Amer.\ Math.\ Soc.\ {\bf 363},
No.\ 1 (2011), 229--257.

\bibitem{MMN-powers} J.\ Migliore, R.\ Mir\'o-Roig and U.\ Nagel, {\em
On the Weak Lefschetz Property for Powers of Linear Forms}, Algebra
Number Theory (to appear).

\bibitem{MN3} J.\ Migliore and U.\ Nagel, {\em Reduced arithmetically
Gorenstein schemes and simplicial polytopes with maximal Betti numbers}, Adv.\ Math.\ {\bf 180} (2003), 1--63.

\bibitem{MN-gor-quad} J.\ Migliore and U.\ Nagel, {\em Gorenstein algebras presented by quadrics}, Preprint, 2011.

\bibitem{MNZ} J.\ Migliore, U.\ Nagel and F.\ Zanello, {\em A characterization of Gorenstein Hilbert functions in codimension four with small initial degree}, Math.\ Res.\ Lett.\ {\bf 15} (2008), no.\ 2, 331--349.

\bibitem{MZ} J.\ Migliore and F.\ Zanello, {\em The Hilbert functions which force the Weak Lefschetz Property}, J.\ Pure Appl.\ Algebra {\bf 210} (2007), 465--471.

\bibitem{MZ2} J.\ Migliore and F.\ Zanello, {\em The strength of the Weak Lefschetz Property}, Illinois J.\ Math.\ {\bf 52}, no.\ 4 (2008), 1417--1433.

\bibitem{NS-Bu}
I.\ Novik, E.\ Swartz, {\em Gorenstein rings through face rings of
manifolds}, Compos.\ Math.\ {\bf 145} (2009),  993–-1000.

\bibitem{PR} K.\ Pardue and B.\ Richert, {\em Syzygies of semi-regular sequences}, Illinois J.\ Math.\ {\bf 53} (2009), 349--364.

\bibitem{RRR} L.\ Reid, L.\ Roberts and M.\ Roitman, {\em On complete
intersections and their Hilbert functions}, Canad.\ Math.\ Bull.\
{\bf 34}  (4) (1991), 525--535.

\bibitem{ss} H.\ Schenck and A.\ Seceleanu, {\em The Weak Lefschetz Property and powers of linear forms in $K[x,y,z]$}, Proc.\ Amer.\ Math.\ Soc.\ {\bf 138} (2010), 2335--2339.

\bibitem{sekiguchi} H.\ Sekiguchi, {\em The upper bound of the Dilworth number and the Rees number of Noetherian local rings  with a Hilbert function}, Adv.\ Math.\ {\bf 124} (1996), 197--206.

\bibitem{ss} S. Seo and H. Srinivasan, {\em On unimodality of Hilbert functions of Gorenstein Artin algebras of embedding dimension four}, preprint 2010.  Available on the arXiv as {\tt arXiv:1011.2732v1}.

\bibitem{stanley-adv} R.\ Stanley,  {\it Hilbert functions
of graded algebras}, Advances in Math.\  {\bf 28} (1978), 57--83.

\bibitem{stanley} R.\ Stanley, {\em Weyl groups, the hard Lefschetz
theorem, and the Sperner property}, SIAM J.\ Algebraic Discrete Methods
{\bf 1} (1980), 168--184.

\bibitem{St-faces} R.\ Stanley, {\em The number of faces of a simplicial
convex polytope}, Advances in  Math.\ {\bf 35} (1980), 236--238.

\bibitem{St-book}
R.\ Stanley, {\em Combinatorics and Commutative Algebra},
    2nd edition. Progress in Mathematics {\bf 41} Birkh\"auser Boston,
     Inc. Boston, MA,
    1996.

\bibitem{watanabe} J.\ Watanabe, {\em The Dilworth number of Artinian
rings and finite posets with rank function}, Commutative Algebra and
Combinatorics, Advanced Studies in Pure Math.\ Vol.\ 11, Kinokuniya Co.\
North Holland, Amsterdam (1987), 303--312.

\bibitem{Z-nonunimodal} F.\ Zanello, {\em A non-unimodal codimension 3 level $h$-vector}, J.\ Algebra {\bf 305} (2006), 949--956.

\bibitem{ZZ} F.\ Zanello and J.\ Zylinski, {\em Forcing the strong Lefschetz and the maximal rank properties}, J.\ Pure Appl.\ Algebra {\bf 213} (2009), 1026--1030.

\end{thebibliography}
\end{document}